\documentclass[11pt, reqno]{amsart}
\usepackage{lmodern}

\usepackage{amsmath,amsfonts,amssymb,amsthm,booktabs,color,epsfig,graphicx,hyperref,url}

\usepackage[table]{xcolor}
\newcommand{\highl}{\cellcolor{gray!25}}

\usepackage{natbib}
\usepackage{dsfont}
\usepackage{upgreek, mathrsfs}
\usepackage{tikz}

\usepackage{rotating}
\usepackage{array} %alignment resp graphics in table
\usepackage{longtable}
\usepackage{enumitem}
\usepackage{calc} 
\allowdisplaybreaks

\usepackage{setspace}
\newcommand{\bbz}{\mathbb{Z}}
\newcommand{\bbn}{\mathbb{N}}

\newcommand{\bfB}{\mathbf{B}}

\newcommand{\bfX}{\mathbf{X}}

\newcommand{\bfY}{\mathbf{Y}}
\newcommand{\bfZ}{\mathbf{Z}}

\newcommand{\bbe}{\mathbb{E}}
\newcommand{\var}{\mathbb{V}}
\newcommand{\bbp}{\mathbb{P}}
\newcommand{\cov}{\text{Cov}}

\newcommand{\bin}{\text{Bin}}
\newcommand{\poi}{\text{Poi}}

\newcommand{\pgf}{\operatorname{pgf}}

\newcommand{\brackets}[1]{\left( #1 \right)}
\newcommand{\bbrackets}[1]{\big( #1 \big)}
\newcommand{\Bbrackets}[1]{\Big( #1 \Big)}

\theoremstyle{plain}% Theorem-like structures provided by amsthm.sty
\newtheorem{theorem}{Theorem}[section]

\newtheorem{proposition}[theorem]{Proposition}

\newtheorem{corollary}[theorem]{Corollary}

\theoremstyle{definition}

\newtheorem{remark}[theorem]{Remark}
\newtheorem{example}[theorem]{Example}

\begin{document}
\title[]{The Integer-valued Moving-Average Random Field} 
\author[A. Silbernagel]{Angelika Silbernagel\textsuperscript{\textasteriskcentered}}
\author[C.H. Wei\ss{}]{Christian H. Wei\ss{}\textsuperscript{\dag}}
%\today

\address{\textsuperscript{\textasteriskcentered} Department of Mathematics and Statistics, Helmut Schmidt University, Hamburg, Germany. \newline
ORCID: \href{https://orcid.org/0009-0002-6993-244X}{\nolinkurl{0009-0002-6993-244X}}.
}
\email{silbernagel@hsu-hh.de}

\address{\textsuperscript{\dag} Department of Mathematics and Statistics, Helmut Schmidt University, Hamburg, Germany. \newline
ORCID: \href{https://orcid.org/0000-0001-8739-6631}{\nolinkurl{0000-0001-8739-6631}}.
}
\email{weissc@hsu-hh.de}

\keywords{Autocorrelation function,
Count data,
Moving-average model,
Random field,
Spatial dependence
}
\subjclass[2020]{Primary 60G60 %Random fields 
62M40, %Random fields; image analysis
Secondary 60G10 %Stationary stochastic processes
62M10 %Time series, auto-correlation, regression, etc. in statistics (GARCH)
}

\begin{abstract}
An integer-valued moving average (INMA) model for count random fields is proposed and investigated. Closed-form expressions are derived for both its marginal distribution and spatial dependence structure, for arbitrary model order and also covering the multilateral case. In particular, general expressions for bivariate distributions and autocovariances are provided. It is shown that the INMA random field can be equipped (among others) with a Poisson marginal distribution. It is also demonstrated that different and well-interpretable dependence structures are possible. For illustration, we discuss a real-world data example and propose an INMA approximation to a given spatial dependence structure. 
\end{abstract}
\maketitle

%\vspace*{-7mm} %%%%%%%%%%%%%%%%%%%%%%%%%%%%%%%%%
\section{Introduction}
\label{section: Introduction}
The first stochastic models for real-valued and continuously distributed (two-dimen\-sional) random fields were developed several decades ago. In particular, ``plane counterparts'' to the classical autoregressive (AR) and moving-average (MA) models known from time series analysis have been discussed extensively by many researchers, such as \citet{whittle54,besag74,haining78,basu93,arezki26}. The modeling of discretely distributed random fields, by contrasts, has been widely ignored for a long time. 
The first integer-valued ARMA-type model for count random fields (i.e., where the outcomes are quantitative with range $\bbn_0=\{0,1,\ldots\}$) has been introduced by \cite{ghodsi12}. There, the authors proposed a spatial analog of the first-order integer-valued AR (INAR$(1)$) model being known from the field of count time series \citep[see, e.g.,][]{mckenzie88,silva05,weiss08,weiss17,moller18,weiss18,altun22,irshad25}. Some more recent achievements on bounded or $\bbz$-valued AR time series models include, e.g., \citet{li25, xu25, zhu26}. More precisely, the count random variables $(X_{s,t}) = (X_{s,t})_{s,t\in\bbz=\{\dots, -1,0,1,\dots\}}$ are located on a two-dimensional regular lattice (grid), and the authors refer to their model as ``first-order spatial INAR (SINAR$(1,1)$)''. However, in what follows, we follow the terminology of \citet{sil_wei_26} and refer to this model as the INAR$(1,1)$ \emph{random field}, in order to avoid confusion with other types of spatial count model and to indicate the grid structure of the generated data.  

\smallskip
The INAR$(1,1)$ random field model of \citet{ghodsi12} has been studied, modified and extended in many subsequent works, including \citet{chutoo21,ghodsi15,sassi23,ghodsi24,karlis24,tabandeh24,ghodsi25,yang25, sil_wei_26}. Interestingly, it seems that the literature solely focused on AR-type models for count random fields so far. For simulating spatially dependent count random fields, however, it would be more convenient to use an MA-type model for count random fields as burn-ins for INAR random fields are computationally expensive, and the simulation of multilateral INAR random fields as considered by \citet{karlis24} is even more demanding. Furthermore, it is generally difficult to derive the stationary marginal distribution of INAR random fields. Both issues could be solved by using integer MA-type random fields instead. To our knowledge, however, such INMA random field models have not been investigated in the literature so far, we are only aware of a few applications within simulation studies, such as in \citet{wei_kim_24}.

\smallskip
Therefore, in what follows, we close this gap in the research literature by adapting the INMA time series model (see \citealp{weiss08} for a detailed discussion) to a unilateral random field structure, and by comprehensively inspecting its stochastic properties. Similar to the time series model, the INMA random field is based on the \emph{binomial thinning operator} denoted by ``$\circ$'', which proved to be an adequate substitute to scalar multiplication for integer-valued time series models \citep{weiss18}. Let $X$ be a discrete count random variable with range $\bbn_0$ and $\alpha \in (0,1)$ a constant. The thinned random variable $\alpha \circ X$ is defined by $\alpha \circ X := \sum^X_{i=1} Z_i$, where $(Z_i)_{i\in\bbn=\{1,2,\ldots\}}$ is a sequence of independent and identically distributed (i.i.d.)\ Bernoulli random variables with $\bbp(Z_i=1)=\alpha$, being independent of $X$. Hence, the conditional distribution of $\alpha \circ X|X$ is the binomial distribution $\bin(X,\alpha)$, what explains the term ``binomial thinning''. Here, the boundary values $\alpha\in\{0,1\}$ can be included by the conventions $0 \circ X := 0$ and $1 \circ X:=X$.

\smallskip
We define the (unilateral) \emph{INMA random field} of order $(q_1,q_2)\in\bbn_0^2$ with $q_1+q_2\geq 1$, abbreviated as INMA$(q_1,q_2)$, by 
\begin{equation}
    \label{eq: model inma}
    X_{s,t} = \sum^{q_1}_{i=0} \sum^{q_2}_{j=0} \beta_{i,j}\, \circ_{s,t}\, \varepsilon_{s-i,t-j},
\end{equation}
where $\beta_{i,j}\in [0,1]$ denote the model parameters, $(\varepsilon_{s,t})=(\varepsilon_{s,t})_{s,t\in\bbz}$ is a sequence of i.i.d.\ non-negative integer-valued random variables with mean $\mu_\varepsilon=\bbe (\varepsilon_{s,t})$ and variance $\sigma^2_\varepsilon=\var(\varepsilon_{s,t})$, and all thinnings concerning \emph{different}~$\varepsilon_{s,t}$ are performed independently of each other and of~$\varepsilon_{s,t}$. 
In \eqref{eq: model inma}, we have added an index ``$s,t$'' below the operator ``$\circ$'' in order to indicate that the thinning is executed at point $(s,t)$; we omit the index whenever it is clear. Note that we allow for the possible dependence of those thinnings that are applied to the \emph{same}~$\varepsilon_{s,t}$, i.e., for the thinnings being summarized in the $(q_1+1)(q_2+1)$-dimensional vectors
\begin{equation}
    \label{eq: INMA representation Y}
    \bfY_{s,t} = (\beta_{0,0}  \circ_{s,t} \varepsilon_{s,t},\ \beta_{0,1}  \circ_{s,t+1} \varepsilon_{s,t},\ \dots,\ \beta_{q_1, q_2} \circ_{s+q_1,t+q_2} \varepsilon_{s,t})^\top.
\end{equation}
More precisely, $\bfY_{s,t}|\varepsilon_{s,t}$ might exhibit (conditional) cross-dependence while different vectors~$\bfY_{s,t}$ are independent of each other (more precisely, the~$(\bfY_{s,t})$ are i.i.d.).
The possible cross-dependencies within~$\bfY_{s,t}$ as well as their consequences on the stochastic properties of the INMA$(q_1,q_2)$ model \eqref{eq: model inma} are further detailed in Section~\ref{section: Special Cases of INMA Random Fields} below. Furthermore, note that we do not need to impose any further independence assumptions between the thinnings and the ``past values'' of $X_{s,t}$ (unlike in the refined INAR random field model according to \citealp{sil_wei_26}), as these are fully implied by the INMA$(q_1,q_2)$ structure.

\smallskip
The remainder of this article is organized as follows. In Sections~\ref{section: Marginal Properties of INMA Random Fields} and~\ref{section: Spatial Dependence in INMA Random Fields}, we investigate fundamental stochastic properties of the INMA random field according to \eqref{eq: model inma}, namely its marginal distribution as well as its spatial dependence structure. In particular, we derive closed-form expressions for the bivariate distributions and autocovariances at a given spatial lag $(k,l)$. In Section~\ref{section: Special Cases of INMA Random Fields}, we present several special cases of the INMA random field, which appear to be particularly relevant for applications. They cover the Poisson INMA random field on the one hand, and special cross-dependence structures for the thinnings in \eqref{eq: INMA representation Y} on the other hand. Section~\ref{section: Numerical Illustrations} presents various illustrations of our novel INMA model for count random fields. We apply the model to a real-world data example on yeast cell counts obtained by a haemacytometer, we use it for approximating the autocorrelation function (ACF) of an INAR random field, and we briefly demonstrate the effect of a multilateral INMA structure. Finally, we conclude in Section~\ref{section: Conclusions and Future Research} and outline possible directions for future research.

%%%%%%%%%%%%%%%%%%%%%%%%
\section{Marginal Properties of INMA Random Fields}
\label{section: Marginal Properties of INMA Random Fields}
In this section, we focus on the marginal distribution of the (strictly stationary) INMA random field according to \eqref{eq: model inma}. It shall get clear that the marginal distribution is fully known and can be expressed in terms of the probability generating function (pgf). In particular, it shall turn out that this marginal distribution is not affected by the actual cross-dependence structure according to \eqref{eq: INMA representation Y}.
The following Proposition~\ref{prop: inma marginal properties} summarizes all relevant marginal properties of $(X_{s,t})$, which are derived in close analogy to the time-series case discussed by \citet[Theorem~2.7]{weiss08}.

\begin{proposition}
    \label{prop: inma marginal properties}
    Let $(X_{s,t})$ be a (strictly) stationary INMA$(q_1,q_2)$ random field satisfying \eqref{eq: model inma}. Then,
    \begin{equation*}
        \mu_X = \bbe (X_{s,t}) = \mu_\varepsilon\, \beta_\bullet \qquad \text{ and } \qquad \sigma_X^2 = \var(X_{s,t}) = \mu_X + (\sigma^2_\varepsilon - \mu_\varepsilon) \sum^{q_1}_{i=0} \sum^{q_2}_{j=0} \beta^2_{i,j},
    \end{equation*}
    where $\beta_\bullet = \sum^{q_1}_{i=0} \sum^{q_2}_{j=0} \beta_{i,j}$. Furthermore, the pgf of $X_{s,t}$, i.e., $\pgf_X(u) = \bbe (u^{X_{s,t})}$, satisfies
    \[
        \pgf_X(u) = \prod^{q_1}_{i=0} \prod^{q_2}_{j=0} \pgf_\varepsilon\bbrackets{1+\beta_{i,j}(u-1)},
    \]
    where $\pgf_\varepsilon$ denotes the pgf of the i.i.d.\ innovations $(\varepsilon_{s,t})$.
\end{proposition}
The proof of Proposition~\ref{prop: inma marginal properties} is provided in Appendix~\ref{app: Proof of Proposition prop: inma marginal properties}. 
Proposition~\ref{prop: inma marginal properties} implies that $X_{s,t}$ is over-/equi-/underdispersed iff $\varepsilon_{s,t}$ is over-/equi-/underdispersed, i.e., the dispersion behavior of $X_{s,t}$ is completely determined by the innovations' marginal distribution.
Moreover, Proposition~\ref{prop: inma marginal properties} shows that having specified the innovations' distribution, we immediately obtain the resulting observations' distribution. For example, see the details in Example~\ref{example: inma poisson} below, Poisson-distributed innovations lead to Poisson observations. 
In particular, inserting $u=0$ into the pgf shows that the probability of observing a zero is given by
\[
    \bbp(X_{s,t}=0) = \pgf_X(0) = \prod^{q_1}_{i=0} \prod^{q_2}_{j=0} \pgf_\varepsilon(1-\beta_{i,j}),
\]
meaning that the zero probability is completely determined by the marginal distribution of the innovations and the choice of parameters $\beta_{i,j}$.
The availability of a closed-form expression for the marginal distribution constitutes a major advantage over the existing INAR random fields, where generally only mean and variance of the observations are known.

\begin{remark}
\label{remark: inma simulation}
Another advantage of INMA compared to INAR random fields is their ease of simulation. If an $n_1\times n_2$ rectangular grid of INMA counts needs to be simulated, $n_1,n_2\in\bbn$, i.e., the counts $X_{1,1}, \dots, X_{n_1,n_2}$, then 
\begin{enumerate}
    \item we first simulate the $(n_1+q_1)\times (n_2+q_2)$ rectangular grid of i.i.d.\ count innovations $\varepsilon_{1-q_1,1-q_2}, \ldots, \varepsilon_{n_1,n_2}$. 
    
    \item Second, for each innovation~$\varepsilon_{s,t}$, we generate the vector~$\bfY_{s,t}$ of all thinnings applied to~$\varepsilon_{s,t}$, recall \eqref{eq: INMA representation Y}.
    
    \item Finally, we create the observations $X_{1,1}, \dots, X_{n_1,n_2}$ by simply adding the appropriate thinnings from the vectors~$(\bfY_{s,t})$, namely according to
\begin{equation}
    \label{eq: model inma Y}
    X_{s,t} = \sum^{q_1}_{i=0} \sum^{q_2}_{j=0} \beta_{i,j}\, \circ_{s,t}\, \varepsilon_{s-i,t-j}
    = \sum^{q_1}_{i=0} \sum^{q_2}_{j=0} Y_{s-i,t-j}^{(i,j)},
\end{equation}
where $Y_{s-i,t-j}^{(i,j)}$ denotes the $(i,j)$\textsuperscript{th} entry of~$\bfY_{s-i,t-j}$.
\end{enumerate}
So only a few additional innovations have to be simulated, but a costly burn-in can be avoided. This differs from the INAR case, where, e.g., \citet{wei_kim_24} used a burn-in period of ``width''~100: they simulated $(n_1+100)\times (n_2+100)$ counts in order to use the $n_1\times n_2$ ``most recent'' counts for their power analyses.

\smallskip
In addition, the sketched simulation approach would be easily extended to a \emph{multilateral INMA random field} of the form
\begin{equation}
    \label{eq: model inma multilateral}
    X_{s,t} = \sum^{q_1}_{i=-p_1} \sum^{q_2}_{j=-p_2} \beta_{i,j}\, \circ_{s,t}\, \varepsilon_{s-i,t-j},
\end{equation}
where the observation~$X_{s,t}$ does not only depend on the ``current'' and ``past'' innovations~$\varepsilon_{s-i,t-j}$ with $i,j\geq 0$, but also on ``future'' innovations~$\varepsilon_{s+k,t+l}$ with $k>0$ or $l>0$, see Figure~\ref{fig: multilateral structure} for illustration. Here, one just has to simulate additional innovations in the first step, the remaining steps are as before. The simulation of the multilateral INAR random field, by contrast, is much more demanding, see \citet[p.~11]{karlis24} for details. Although the present research mainly focuses on the unilateral INMA model according to \eqref{eq: model inma} (formally, the multilateral model \eqref{eq: model inma multilateral} might be understood as a ``shifted'' unilateral INMA$(p_1+q_1,p_2+q_2)$ model), we shall sometimes also briefly refer to the above multilateral extension, see Remark~\ref{remark: multilateral} and Section~\ref{section: Numerical Illustrations}.

\begin{figure}
    \centering
    a)\scalebox{1}{
   \begin{tikzpicture}
        % Unilateral
        \node (p1) at (0,2) {$\varepsilon_{s-1,t-1}$};
        \node (p2) at (0,3) {$\varepsilon_{s-1,t}$};
        \node (p4) at (2,2) {$\varepsilon_{s,t-1}$};
        \node (p5) at (2,3) {$X_{s,t}$};
        \draw[thick, ->] (p1.north east) -- (p5.south west);
        \draw[thick, ->] (p2.east) -- (p5.west);
        \draw[thick, ->] (p4.north) -- (p5.south);
        \node at (1,1.2) {unilateral};
    \end{tikzpicture}}
    \qquad
    b)\scalebox{1}{
    \begin{tikzpicture}
        % Queen
        \node (p1) at (0,2) {$\varepsilon_{s-1,t-1}$};
        \node (p2) at (0,3) {$\varepsilon_{s-1,t}$};
        \node (p3) at (0,4) {$\varepsilon_{s-1,t+1}$};
        \node (p4) at (2,2) {$\varepsilon_{s,t-1}$};
        \node (p5) at (2,3) {$X_{s,t}$};
        \node (p6) at (2,4) {$\varepsilon_{s,t+1}$};
        \node (p7) at (4,2) {$\varepsilon_{s+1,t-1}$};
        \node (p8) at (4,3) {$\varepsilon_{s+1,t}$};
        \node (p9) at (4,4) {$\varepsilon_{s+1,t+1}$};
        \draw[thick, ->] (p1.north east) -- (p5.south west);
        \draw[thick, ->] (p2.east) -- (p5.west);
        \draw[thick, ->] (p3.south east) -- (p5.north west);
        \draw[thick, ->] (p4.north) -- (p5.south);
        \draw[thick, ->] (p6.south) -- (p5.north);
        \draw[thick, ->] (p7.north west) -- (p5.south east);
        \draw[thick, ->] (p8.west) -- (p5.east);
        \draw[thick, ->] (p9.south west) -- (p5.north east);
        \node at (2,1.2) {multilateral};
    \end{tikzpicture}
    }
    \caption{Comparison of a) unilateral to b) multilateral dependence structure in case of first-order INMA random field.}
    \label{fig: multilateral structure}
\end{figure}

\smallskip
So altogether, the INMA random fields offer a simple way of simulating spatial dependence in regular grid data, also see Section~\ref{section: Numerical Illustrations} for a further discussion. Hence, they are useful for, e.g., simulation-based power analyses of tests for spatial dependence, such as in \citet{wei_kim_24}.
\end{remark}

%%%%%%%%%%%%%%%%%%%%%%%%%%%
\section{Spatial Dependence in INMA Random Fields}
\label{section: Spatial Dependence in INMA Random Fields}
Already the first-order INMA$(1,1)$ random field
\[
    X_{s,t} = \beta_{0,0}\circ\varepsilon_{s,t} + \beta_{1,0}\circ\varepsilon_{s-1,t} + \beta_{0,1}\circ\varepsilon_{s,t-1} + \beta_{1,1}\circ\varepsilon_{s-1,t-1}
\]
involves more than one thinning operator at a location~$(s,t)$, even if we set $\beta_{0,0}\equiv1$. Consequently, each innovation $\varepsilon_{s,t}$ is involved in the four thinnings $\beta_{0,0}\circ_{s,t}\varepsilon_{s,t}$, $\beta_{1,0}\circ_{s+1,t}\varepsilon_{s,t}$, $\beta_{0,1}\circ_{s,t+1}\varepsilon_{s,t}$, and $\beta_{1,1}\circ_{s+1,t+1}\varepsilon_{s,t}$. For the general INMA$(q_1,q_2)$ random field, $\varepsilon_{s,t}$ is even involved in $(q_1+1)(q_2+1)$ thinnings, namely $\beta_{k,l}\circ_{s+k,t+l}\varepsilon_{s,t}$ with $0\leq k \leq q_1$, $0 \leq l \leq q_2$ as summarized within the vector~$\bfY_{s,t}$ from \eqref{eq: INMA representation Y}. These thinning operators are probabilistic, so in analogy to the classical INMA$(q)$ time series model, their joint conditional distribution $(\beta_{k,l}\circ_{s+k,t+l}\varepsilon_{s,t}, 0\leq k \leq q_1, 0 \leq l \leq q_2\, |\, \varepsilon_{s,t})$ has to be specified, which may lead to different types of model for the same model order $(q_1,q_2)$, see \citet{weiss08} for a comprehensive discussion in the time series case. Altogether, the final spatial dependence structure of an INMA$(q_1,q_2)$ random field is not only determined by the actual MA-structure, but also by the possible cross-dependencies within the vectors~$\bfY_{s,t}$.

\smallskip
Let us first take a bird’s-eye view of the INMA$(q_1,q_2)$ model’s dependency structure, namely by considering a full $n_1\times n_2$ rectangular grid of observations, $n_1,n_2\in\bbn$, which we may summarize in the $(n_1+q_1) (n_2+q_2)$-dimensional vector $\bfX=(X_{1,1}, \dots, X_{n_1,n_2})^\top$. As already explained in Remark~\ref{remark: inma simulation}, these observations are determined once 
\begin{enumerate}
    \item the $(n_1+q_1)\times (n_2+q_2)$ i.i.d.\ count innovations $\varepsilon_{1-q_1,1-q_2}, \ldots, \varepsilon_{n_1,n_2}$ have been generated, and
    \item all corresponding thinnings have been executed. 
\end{enumerate}
In other words, as soon as the $(n_1+q_1)\times (n_2+q_2)$ vectors~$\bfY_{1-q_1,1-q_2}, \ldots, \bfY_{n_1,n_2}$, each being of dimension $(q_1+1)(q_2+1)$, have been generated, the counts in~$\bfX$ are determined via \eqref{eq: model inma Y}. This relation between~$\bfX$ and the~$\bfY_{s,t}$ can be expressed by a simple matrix multiplication. Let $\mathbf{Y}=(\bfY_{1-q_1,1-q_2}, \dots, \bfY_{n_1,n_2})^\top$ be the $(q_1+1)(q_2+1)(n_1+q_1)(n_2+q_2)$-dimensional vector being obtained by merging all individual vectors~$\bfY_{s,t}$. Then, \eqref{eq: model inma Y} implies an $(n_1 n_2 \times (q_1+1)(q_2+1)(q_1+n_1)(q_2+n_2))$-dimensional matrix $\bfB$ consisting of only the values $0$ and $1$ such that
\begin{equation}
    \label{eq: INMA representation}
    \bfX = \bfB \cdot \bfY.
\end{equation}
For illustration purposes, we state the representation for the case $n_1=n_2=2$ and $q_1=q_2=1$. Then, \eqref{eq: INMA representation} simplifies to
\[
    \begin{pmatrix}
        X_{1,1} & X_{1,2} & X_{2,1} & X_{2,2}
    \end{pmatrix}^\top
    \ =\ 
    \bfB \cdot 
    \begin{pmatrix}
        \bfY_{0,0} & \bfY_{0,1} & \bfY_{0,2} & \bfY_{1,0} & \bfY_{1,1} & \bfY_{1,2} & \bfY_{2,0} & \bfY_{2,1} & \bfY_{2,2}
    \end{pmatrix}^\top
\]
with $\bfY_{s,t} = (\beta_{0,0}  \circ_{s,t} \varepsilon_{s,t},\ \beta_{1,0}  \circ_{s+1,t} \varepsilon_{s,t},\ \beta_{0,1}  \circ_{s,t+1} \varepsilon_{s,t},\ \beta_{1,1} \circ_{s+1,t+1} \varepsilon_{s,t})^\top$ and
\[
\bfB = \begin{pmatrix}
    0001 & 0100 & 0000 & 0010 & 1000 & 0000 & 0000 & 0000 & 0000 \\
    0000 & 0001 & 0100 & 0000 & 0010 & 1000 & 0000 & 0000 & 0000 \\
    0000 & 0000 & 0000 & 0001 & 0100 & 0000 & 0010 & 1000 & 0000 \\
    0000 & 0000 & 0000 & 0000 & 0001 & 0100 & 0000 & 0010 & 1000 
\end{pmatrix}.
\]
We can see, for example, that the components of $\bfY_{1,1}$ (corresponding to the fifth block of columns in~$\bfB$) appear in all four observations, which may lead to dependence between these observations according to the actual model specifications. It is also clear that the distribution of~$\bfX = \bfB\cdot\bfY$ is fully determined by the distribution of~$\bfY$, where the latter can be given after having fixed some notation. Let $(Z_{s,t;r}^{(i,j)})_{1\leq r \leq \varepsilon_{s,t}}$ denote the counting series involved in the thinning applied to $\varepsilon_{s,t}$ at point $(s+i,t+j)$ for $0 \leq i \leq q_1$ and $0 \leq j \leq q_2$, so $\beta_{i,j}\circ_{s+i,t+j} \varepsilon_{s,t} = \sum_{r=1}^{\varepsilon_{s,t}} Z_{s,t;r}^{(i,j)}$ with $\bbp(Z_{s,t;r}^{(i,j)}=1)=\beta_{i,j}$. Then, the $(q_1+1)(q_2+1)$-dimensional vectors $\bfZ_{s,t;r}:=(Z_{s,t;r}^{(0,0)}, \ldots, Z_{s,t;r}^{(q_1,q_2)})^\top$ are i.i.d.\ for all $s,t$ and $r$, but there might certainly be dependence between the components within~$\bfZ_{s,t;r}$. Now, $\bfY$ consists of mutually i.i.d.\ vectors~$\bfY_{s,t}$, where $\bfY_{s,t} = \sum_{r=1}^{\varepsilon_{s,t}} \bfZ_{s,t;r}$. So the distribution of~$\bfY$ follows from the pgf of~$\bfY_{s,t}$, which is given by
\begin{eqnarray}
\nonumber
\pgf_{\bfY_{s,t}}(u_{00},\ldots,u_{q_1q_2}) &=& \bbe\Big(\cdots u_{ij}^{Y_{s,t;r}^{(i,j)}}\cdots\Big)
= \bbe\Big[\bbe\Big(\cdots u_{ij}^{\sum_{r=1}^{\varepsilon_{s,t}} Z_{s,t;r}^{(i,j)}}\cdots \Big| \varepsilon_{s,t}\Big)\Big]
\\
\nonumber
&=& \bbe\Big[\bbe\Big(\cdots u_{ij}^{Z_{s,t;1}^{(i,j)}}\cdots \Big)^{\varepsilon_{s,t}}\Big]
\\
\label{equation: pgf Y}
&=& \pgf_{\varepsilon}\big(\pgf_{\bfZ}(u_{00},\ldots,u_{q_1q_2})\big),
\end{eqnarray}
where $\pgf_{\bfZ}$ refers to the pgf of $\bfZ_{s,t;r}$.
Let us now take a closer look at the INMA$(q_1,q_2)$ model’s dependency structure. 
Following the approach of \citet{weiss08}, we derive closed-form expressions for the autocovariance function (ACvF) and bivariate pgf without specifying the cross-dependence within~$\bfY_{s,t}$ or~$\bfZ_{s,t;r}$, respectively. Only later in Section~\ref{section: Special Cases of INMA Random Fields}, we consider possible choices of distributions for~$\bfZ_{s,t;r}$, which are natural in the sense that they constitute natural generalizations of established models in the time series case. Our main theorem is as follows.

\begin{theorem}
    \label{theorem: inma dependence}
    Let $(X_{s,t})$ be a stationary INMA$(q_1,q_2)$ random field satisfying \eqref{eq: model inma}, and define the index sets
    $$
    \mathcal{S}_{kl}\ :=\ \big\{(i+k,j+l) :\quad 0 \leq i \leq q_1,\ 0 \leq j \leq q_2\big\}, \qquad k,l\in\bbz.
    $$
    \begin{itemize}
        \item[(i)] The bivariate pgf with spatial lag $(k,l)$, $\pgf(u_1,u_2;k,l) := \bbe(u_1^{X_{s,t}} u_2^{X_{s-k,t-l}})$, is given by the expression
    \begin{align*}
        \pgf(u_1,u_2;k,l) = & \prod_{(i,j) \in \mathcal{S}_{00}\setminus(\mathcal{S}_{00}\cap \mathcal{S}_{kl})} \pgf_\varepsilon\Big(1+\beta_{i,j}\,(u_1-1)\Big) \\
        &\ \cdot \prod_{(i,j) \in \mathcal{S}_{kl}\setminus(\mathcal{S}_{00}\cap \mathcal{S}_{kl})}  \pgf_\varepsilon\Big(1+\beta_{i-k,j-l}\,(u_2-1)\Big) 
        \\
        &\ \cdot \prod_{(i,j) \in \mathcal{S}_{00}\cap \mathcal{S}_{kl}}  \pgf_\varepsilon\Big(1+\beta_{i,j}\,(u_1-1) + \beta_{i-k,j-l}\,(u_2-1) \\
        &\hspace{15ex} + \bbp(Z^{(i,j)}_{s,t;1}=Z^{(i-k,j-l)}_{s,t;1}=1) \cdot(u_1-1)(u_2-1)\Big),
    \end{align*}
    which covers the marginal pgf from Proposition~\ref{prop: inma marginal properties} (with $u=u_1 u_2$) for $k=l=0$ (using the convention that empty products are equal to~1).

    \item[(ii)] The ACvF satisfies
    \begin{align}
        \gamma(k,l) = (\sigma^2_\varepsilon &- \mu_\varepsilon) \sum_{(i,j)\in\mathcal{S}_{00}\cap\mathcal{S}_{kl}} \beta_{i,j}\beta_{i-k,j-l} \nonumber \\
        &+ \mu_\varepsilon \sum_{(i,j)\in\mathcal{S}_{00}\cap\mathcal{S}_{kl}} \bbp(Z_{s-i,t-j;1}^{(i,j)}=Z_{s-i,t-j;1}^{(i-k,j-l)}=1), \label{eq: inma sacf}
    \end{align}
    which covers the variance formula from Proposition~\ref{prop: inma marginal properties} for $k=l=0$ (using the convention that empty sums are equal to $0$).
    \end{itemize}
\end{theorem}
Note that Theorem~1 does not require additional assumptions apart from the model definition (including the assumptions on the innovations). The proof of Theorem~\ref{theorem: inma dependence} is provided in Appendix~\ref{app: Proof of Theorem theorem: inma dependence}. Note the analogy of Theorem~\ref{theorem: inma dependence} to Theorems~2.11 and~2.12 in \citet{weiss08}, i.e., the INMA dependence structure directly extends from the time-series to the random-fields case. The latter is illustrated by Figure~\ref{fig: proof inma bivariate pgf} which depicts the innovations $\varepsilon_{s-i,t-j}$ appearing in $X_{s,t}$ (dashed rectangle) and $X_{s-k,t-l}$ (dotted rectangle), respectively. The overlap between the respective innovations, which corresponds to spatial dependence between $X_{s,t}$ and $X_{s-k,t-l}$, is highlighted in gray. This overlap shows that spatial dependence only exists if $|k|\leq q_1$ and $|l|\leq q_2$. If one of these conditions is violated, there is no gray overlap region anymore, implying the independence of~$X_{s,t}$ and~$X_{s-k,t-l}$. 

\begin{figure}
    \centering
    a)%
    \scalebox{0.7}{
    \begin{tikzpicture}
        \draw[->, very thick, gray] (0.1,0.1) -- (-6,0.1);
        \draw[->, very thick, gray] (0.1,0.1) -- (0.1,-5);
        \draw[thick, dashed] (0,0) -- (-4,0) -- (-4,-3) -- (0,-3) -- (0,0);
        \draw[thick, dotted] (-1.5,-1) -- (-5.5,-1) -- (-5.5,-4) -- (-1.5,-4) -- (-1.5,-1);
        \fill[gray!30] (-1.5,-1) rectangle (-4,-3);

        \draw (-1.5,0) -- (-1.5,0.2);
        \draw (-4,0) -- (-4,0.2);
        \draw (-5.5,0) -- (-5.5,0.2);
        \node[anchor=south] at (-1.5,0.2) {$s-k$};
        \node[anchor=south] at (-4,0.2) {$s-q_1$};
        \node[anchor=south] at (-5.5,0.2) {$s-k-q_1$};

        \draw (0,-1) -- (0.2,-1);
        \draw (0,-3) -- (0.2,-3);
        \draw (0,-4) -- (0.2,-4);
        \node[anchor=west] at (0.2,-1) {$t-l$};
        \node[anchor=west] at (0.2,-3) {$t-q_2$};
        \node[anchor=west] at (0.2,-4) {$t-l-q_2$};

        \filldraw[black] (0,0) circle (2pt) node[anchor=west]{\ $(s,t)$};
        \filldraw[black] (-4,-3) circle (2pt) node[anchor=south west]{$(s-q_1,t-q_2)$};
        \filldraw[black] (-1.5,-1) circle (2pt) node[anchor=south]{$(s-k,t-l)$};

        \node[anchor=north] at (-3,-0.2) {$\mathcal{S}_{00}$};
        \node[anchor=south] at (-4.5,-3.8) {$\mathcal{S}_{kl}$};
        \node at (-2.75,-2) {$\mathcal{S}_{00}\cap \mathcal{S}_{kl}$};
    \end{tikzpicture}}\qquad
    b)%
    \scalebox{0.7}{
    \begin{tikzpicture}
        \draw[<->, very thick, gray] (2.5,0.1) -- (-4.5,0.1);
        \draw[->, very thick, gray] (0.1,0.1) -- (0.1,-5);
        \draw[thick, dashed] (0,0) -- (-4,0) -- (-4,-3) -- (0,-3) -- (0,0);
        \draw[thick, dotted] (2,-1) -- (-2,-1) -- (-2,-4) -- (2,-4) -- (2,-1);
        \fill[gray!30] (0,-1) rectangle (-2,-3);

        \draw (2,0) -- (2,0.2);
        \draw (-4,0) -- (-4,0.2);
        \draw (-2,0) -- (-2,0.2);
        \node[anchor=south] at (2,0.2) {$s-k$};
        \node[anchor=south] at (-4,0.2) {$s-q_1$};
        \node[anchor=south] at (-2,0.2) {$s-k-q_1$};

        \draw (0,-1) -- (0.2,-1);
        \draw (0,-3) -- (0.2,-3);
        \draw (0,-4) -- (0.2,-4);
        \node[anchor=south west] at (0.2,-1) {$t-l$};
        \node[anchor=west] at (0.2,-3) {$t-q_2$};
        \node[anchor=south west] at (0.2,-4) {$t-l-q_2$};

        \filldraw[black] (0,-0.1) circle (2pt) node[anchor=west]{\ $(s,t)$};
        \filldraw[black] (0,-1) circle (2pt) node[anchor=south east]{$(s,t-l)$};
        \filldraw[black] (-2,-3) circle (2pt) node[anchor=north east]{$(s-k-q_1,t-q_2)$};

        \node[anchor=north] at (-3,-0.7) {$\mathcal{S}_{00}$};
        \node[anchor=north] at (1,-1.7) {$\mathcal{S}_{kl}$};
        \node at (-1,-2) {$\mathcal{S}_{00}\cap \mathcal{S}_{kl}$};
    \end{tikzpicture}}
    \caption{Illustration of the regions for the innovations $\varepsilon_{s-i,t-j}$ appearing in $X_{s,t}$ (with $(i,j)\in\mathcal{S}_{00}$, dashed) and $X_{s-k,t-l}$ (with $(i,j)\in\mathcal{S}_{kl}$, dotted) for a) $k,l\geq 0$ and b) $k <0$ and $l\geq 0$. The region referring to the intersection $\mathcal{S}_{00}\cap\mathcal{S}_{kl}$ is highlighted in gray.}
    \label{fig: proof inma bivariate pgf}
\end{figure}

\begin{remark}
\label{remark: multilateral}
    It is easy to check that the properties from Proposition~\ref{prop: inma marginal properties} and Theorem~\ref{theorem: inma dependence} can be generalized to the multilateral structure presented in \eqref{eq: model inma multilateral} as the method of proof carries over. In this case, the sums and products appearing in Proposition~\ref{prop: inma marginal properties} have to be adjusted to $-p_1\leq i \leq q_1$ and $-p_2 \leq j \leq q_2$, including the sum appearing in the definition of $\beta_\bullet$. The assertions from Theorem~\ref{theorem: inma dependence} remain valid if $\mathcal{S}_{kl}$ is replaced with 
    $$
        \mathcal{S}^\ast_{kl}\ :=\ \big\{(i+k,j+l) :\quad -p_1 \leq i \leq q_1,\ -p_2 \leq j \leq q_2\big\}
    $$
    for all $k,l \in \bbz$. 
\end{remark}

%%%%%%%%%%%%%%%%%%%%%%%%
\section{Special Cases of INMA Random Fields}
\label{section: Special Cases of INMA Random Fields}

In this section, we present a couple of important special cases. First, we consider particular types of marginal distribution. A natural distributional choice for INMA models is the compound Poisson (CP) family \citep[Section~2.1.3]{weiss18}, which covers the popular Poisson and negative-binomial distribution as special cases. Due to the CP's invariance properties given in Lemma~2.1.3.2 of \citet{weiss18}, namely the invariance with respect to addition and binomial thinning, it is clear that CP-distributed innovations~$\varepsilon_{s,t}$ in \eqref{eq: model inma} also imply CP-distributed observations~$X_{s,t}$. Here, results become particularly simple in the Poisson case.

\begin{example}
    \label{example: inma poisson}
    We refer to \eqref{eq: model inma} as the \emph{Poisson INMA$(q_1,q_2)$ random field} if its innovations are Poisson-distributed, $\varepsilon_{s,t}\sim\text{Poi}(\mu_\varepsilon)$ with $\mu_\varepsilon>0$. Then, the pgf of $\varepsilon_{s,t}$ is given by \citep[Example~A.1.1]{weiss18}
    \begin{equation}
        \label{eq: pgf poisson}
        \pgf_\varepsilon(u) 
        %= e^{-\mu_\varepsilon} \sum^\infty_{k=0} \frac{\mu_\varepsilon^k}{k!} u^k 
        = \exp\big(\mu_\varepsilon(u-1)\big).
    \end{equation}
    Thus, it follows from Proposition~\ref{prop: inma marginal properties} that
    \begin{align*}
        \pgf_X(u) &= \prod^{q_1}_{i=0} \prod^{q_2}_{j=0} \pgf_\varepsilon(1+\beta_{i,j}(u-1)) \\
        &= \prod^{q_1}_{i=0} \prod^{q_2}_{j=0} \exp\big(\mu_\varepsilon \beta_{i,j} (u-1)\big) = \exp\big(\mu_\varepsilon \beta_\bullet (u-1)\big).
    \end{align*}
    Since the pgf uniquely encodes the probability mass function (PMF) of a distribution, it follows that $X_{s,t}\sim\text{Poi}(\mu_X)$ with $\mu_X = \mu_\varepsilon \beta_\bullet = \bbe (X_{s,t}) = \var(X_{s,t})$ \citep[Example~A.1.1]{weiss18}.
\end{example}
Example~\ref{example: inma poisson} shows that the Poisson INMA$(q_1,q_2)$ random field has both Poisson-distributed innovations and observations. The possibility of having a Poisson \linebreak marginal distribution is a major advantage over existing INAR-type models for random fields (see Section~\ref{section: Numerical Illustrations} for further discussion), where only some moments of the marginal distribution are known so far. 

\smallskip
Also the Poisson INMA$(q_1,q_2)$'s spatial dependence structure is characterized by simple closed-form expressions, as summarized by the following corollary to Theorem~\ref{theorem: inma dependence}.

\begin{corollary}
    \label{corollary: poisson inma dependence}
    Let $(X_{s,t})$ be a stationary Poisson INMA$(q_1,q_2)$ random field according to Example~\ref{example: inma poisson}. 
    
    \begin{itemize}
        \item[(i)] The ACF from Theorem~\ref{theorem: inma dependence} becomes 
    \begin{equation}
    \label{PoiINMA_acf}
        \rho(k,l) = \frac{1}{\beta_\bullet}
            \sum_{(i,j)\in\mathcal{S}_{00}\cap\mathcal{S}_{kl}} \bbp(Z_{s-i,t-j;1}^{(i,j)}=Z_{s-i,t-j;1}^{(i-k,j-l)}=1).
    \end{equation}
    
        \item[(ii)] The bivariate pgf $\pgf(u_1,u_2;k,l)$ from Theorem~\ref{theorem: inma dependence} simplifies to 
    \[
        \pgf(u_1,u_2;k,l)= \exp\big(\mu_X\, (u_1+u_2-2)\big)\, \exp\big(\mu_X\, \rho(k,l)\, (u_1-1)(u_2-1)\big),
    \]
    which is the pgf of a bivariate Poisson distribution \citep[Example~A.3.1]{weiss18}.
    
        \item[(iii)] Conditional mean and variance are given by
    \begin{eqnarray*}
        \bbe(X_{s,t}\, |\, X_{s-k,t-l}=x) &=& \mu_X\, \big(1-\rho(k,l)\big) + \rho(k,l)\cdot x 
    \\
    \text{and}\quad 
        \var(X_{s,t}\, |\, X_{s-k,t-l}=x) &=& \big(1-\rho(k,l)\big)\, \big(\mu_X +\rho(k,l)\cdot x\big).
    \end{eqnarray*}
        
        \item[(iv)] The distribution of ``jumps'' between~$X_{s,t}$ and~$X_{s-k,t-l}$ is determined by 
    \begin{equation*}
        \bbp(X_{s,t}=X_{s-k,t-l}\pm j) = \exp\Big(-2\mu_\varepsilon\beta_\bullet \big(1-\rho(k,l)\big)\Big) \cdot I_j\Big(2\mu_\varepsilon\beta_\bullet \big(1-\rho(k,l)\big)\Big),
    \end{equation*}
    where $I_j(z) = \sum^\infty_{r=0} (\tfrac{z}{2})^{2r+j}/\big(r! \cdot (r+j)!\big)$, $j \in \bbn_0$,  denotes the \emph{modified Bessel function of the first kind}. Therefore,
    \begin{align*}
        \bbp(X_{s,t}<X_{s-k,t-l}) &= \bbp(X_{s,t}>X_{s-k,t-l}) \\
        &= \frac{1-\exp\Big(-2\mu_\varepsilon\beta_\bullet \big(1-\rho(k,l)\big)\Big) \cdot I_0\Big(2\mu_\varepsilon\beta_\bullet \big(1-\rho(k,l)\big)\Big)}{2}.
    \end{align*}
    \end{itemize}
\end{corollary}
The proof of Corollary~\ref{corollary: poisson inma dependence} is provided in Appendix~\ref{app: Proof of Corollary corollary: poisson inma dependence}. Again, it is worth pointing out the analogy to the time-series case, see Examples~2.8 and~2.13 as well as Theorem~3.3 in \citet{weiss08}. Part~(iii) of Corollary~\ref{corollary: poisson inma dependence} could be useful for model diagnostics, by computing a kind of standardized Pearson residuals \citep[see][Section~2.4]{weiss18}, while part~(iv) could be utilized for alternative dependence measures, see the discussion of Cohen's~$\kappa$ in \citet{weiss08}.

\bigskip
The next group of special cases refers to the possible dependence of those thinnings that are applied to the same innovation~$\varepsilon_{s,t}$ (recall that the marginal properties discussed in Section~\ref{section: Marginal Properties of INMA Random Fields} are not affected by such dependencies). The subsequent discussion is inspired by that of the time-series case in \citet{weiss08}, where the~$\varepsilon_{s,t}$ are understood as expressing the size of some kind of ``population'' at point $(s,t)$. A binomial thinning applied to~$\varepsilon_{s,t}$ is then understood as the outcome of some random Bernoulli experiment being applied to the individuals of the $(s,t)$\textsuperscript{th} population (``generation~$(s,t)$''), namely as the number of individuals being ``successful'' within this experiment. Then, different experimental scenarios lead to different cross-dependencies between the thinnings being applied to~$\varepsilon_{s,t}$, i.e., to different cross-dependencies within~$\bfY_{s,t}$ and~$\bfZ_{s,t;r}$, respectively. In \citet[Appendix~B]{weiss08}, altogether four such scenarios are discussed for the time-series case, where the following two cases extend naturally to random fields.

\begin{example}
    \label{example: inma independence}
    In accordance to the \emph{independence model} for the classical INMA($q$) process in \citet{weiss08}, which was originally proposed by \cite{mckenzie88}, our first special case assumes that the components of $\bfZ_{s,t;r}$ are mutually independent, i.e., the aforementioned experiment is applied to the individuals of generation~$(s,t)$ each location $(s+i,t+j)$ anew, $i=0,\ldots,q_1$ and $j=0,\ldots,q_2$. 
    In this case,
    \begin{align*}
        \bbp(Z_{s-i,t-j;1}^{(i,j)}=Z_{s-i,t-j;1}^{(i-k,j-l)}=1) &= \bbp(Z_{s-i,t-j;1}^{(i,j)}=1) \bbp(Z_{s-i,t-j;1}^{(i-k,j-l)}=1) \\
        &= \beta_{i,j} \beta_{i-k,j-l}
    \end{align*}
    for $(k,l) \neq (0,0)$, and the ACvF in \eqref{eq: inma sacf} simplifies to 
    \[
        \gamma(k,l) = \sigma^2_\varepsilon \sum_{(i,j)\in\mathcal{S}_{00}\cap\mathcal{S}_{kl}} \beta_{i,j}\beta_{i-k,j-l} \qquad \text{for } (k,l)\neq(0,0).
    \]
    This further simplifies to
    \[
        \rho(k,l) = \sum_{(i,j)\in\mathcal{S}_{00}\cap\mathcal{S}_{kl}} \beta_{i,j}\beta_{i-k,j-l}/\beta_\bullet \qquad \text{for } (k,l)\neq(0,0)
    \]
    if the innovations $\varepsilon_{s,t}\sim\text{Poi}(\mu_\varepsilon)$ follow a Poisson distribution, recall \eqref{PoiINMA_acf}.
\end{example}
The second special case is inspired by the so-called \emph{sale model} in \citet{weiss08}, which was originally proposed by \cite{brannas01}. 
Here, the original intuition was to assume the innovations to be the number of new items of some product with a fixed shelf life (determining the order of the moving average), where each item can be sold at most once. 
While the original motivation behind the sale model clearly assumes a time context, we can adapt it to the random-fields case by thinking of regional spread in an agricultural setup. For example, assume that seeds are spread over a field, with the number of seeds sown per grid cell being~$\varepsilon_{s,t}$. Due to wind, some of the seeds spread to neighboring grid cells, which is expressed by the thinnings applied to~$\varepsilon_{s,t}$. Then, $X_{s,t}$ is the total number of seeds finally arriving in grid cell~$(s,t)$.

\begin{example}
    \label{example: inma spread model}
    For the ``spread model'', let $\beta_\bullet\leq1$, and suppose that at most one component in~$\bfZ_{s,t;r}$ is equal to one, all others are equal to zero, and $\bbp(Z_{s,t;r}^{(i,j)}=1) = \beta_{i,j}$. In particular, $\bfZ_{s,t;r}$ becomes the zero vector with probability $1-\beta_\bullet$. This assumption corresponds to each ``individual'' being ``active'' in at most one grid cell, like it happens for seed being spread over a field. For this model, it holds that $\bbp(Z_{s-i,t-j;1}^{(i,j)}=Z_{s-i,t-j;1}^{(i-k,j-l)}=1) =0$ for $(k,l)\neq (0,0)$, since it is impossible for an individual to be active twice. Hence, the ACvF in \eqref{eq: inma sacf} simplifies to
    \[
        \gamma(k,l) = (\sigma^2_\varepsilon - \mu_\varepsilon) \sum_{(i,j)\in\mathcal{S}_{00}\cap\mathcal{S}_{kl}} \beta_{i,j}\beta_{i-k,j-l}
        \quad\text{for }
        (k,l)\neq (0,0).
    \]
    So the ACvF gets reduced by the factor $(\sigma^2_\varepsilon - \mu_\varepsilon)/\sigma^2_\varepsilon\ <1$ compared to the independence model in Example~\ref{example: inma independence}. 
    If the innovations are additionally Poisson distributed (with $\sigma^2_\varepsilon = \mu_\varepsilon$), then parts~(i) and~(ii) of Corollary~\ref{corollary: poisson inma dependence} even imply that spatial dependence vanishes, so the $X_{s,t}\sim \text{Poi}(\mu_X)$ are simply i.i.d.
\end{example}
The remaining two special cases of INMA$(q)$ time series models, the ``lifetime model'' of \citet{brannas01} and the ``changing-states model'' of \citet{alosh88}, are difficult to extend beyond a time context and are, thus, not further considered here.

\smallskip
Finally, let us have a look at the vectors~$\bfY_{s,t}$ of all thinnings being applied to~$\varepsilon_{s,t}$, the pgf of which is given in \eqref{equation: pgf Y}:
$$
\pgf_{\bfY_{s,t}}(u_{00},\ldots,u_{q_1q_2}) = \pgf_{\varepsilon}\big(\pgf_{\bfZ}(u_{00},\ldots,u_{q_1q_2})\big).
$$
This pgf is necessary if one wants to compute the joint probability of a full rectangle $\bfX = \bfB \cdot \bfY$ of INMA$(q_1,q_2)$-counts, recall \eqref{eq: INMA representation}. For the independence model according to Example~\ref{example: inma independence}, we have
\begin{align}
\pgf_{\bfZ}(u_{00},\ldots,u_{q_1q_2}) 
= \bbe\Big(u_{00}^{Z_{s,t;r}^{(0,0)}}\cdots u_{q_1q_2}^{Z_{s,t;r}^{(q_1,q_2)}}\Big) 
&= \bbe\Big(u_{00}^{Z_{s,t;r}^{(0,0)}}\Big)\cdots \bbe\Big(u_{q_1q_2}^{Z_{s,t;r}^{(q_1,q_2)}}\Big) \nonumber \\
&= \prod_{i=0}^{q_1} \prod_{j=0}^{q_2} (1-\beta_{i,j}+\beta_{i,j} u_{ij}), \label{equation: pgf Z independence}
\end{align}
while for the spread model according to Example~\ref{example: inma spread model}, we have
\begin{equation}
\label{equation: pgf Z spread}
\pgf_{\bfZ}(u_{00},\ldots,u_{q_1q_2}) 
= 1-\beta_\bullet + \sum_{i=0}^{q_1} \sum_{j=0}^{q_2} \beta_{i,j} u_{ij}.
\end{equation}
If considering a Poisson INMA$(q_1,q_2)$ random field according to Example~\ref{example: inma poisson} with $\pgf_\varepsilon(u) = \exp\big(\mu_\varepsilon(u-1)\big)$, recall \eqref{eq: pgf poisson}, then \eqref{equation: pgf Z independence} and \eqref{equation: pgf Z spread} further simplify. In particular, for the spread model, we obtain
\begin{align}
\pgf_{\bfY}(u_{00},\ldots,u_{q_1q_2}) 
&= \exp\Big[\mu_\varepsilon\Big(\sum_{i=0}^{q_1} \sum_{j=0}^{q_2} \beta_{i,j} (u_{ij}-1)\Big)\Big] \nonumber \\
&= \prod_{i=0}^{q_1} \prod_{j=0}^{q_2} \exp\big(\mu_\varepsilon \beta_{i,j} (u_{ij}-1)\big), \label{equation: pgf Y spread}
\end{align}
so~$\bfY_{s,t}$ (and thus~$\bfX$) consists of independent Poisson counts.

\section{Numerical Illustrations}
\label{section: Numerical Illustrations}
\begin{figure}[t]
\centering
a)\includegraphics[viewport=30 45 235 240, clip=, scale=0.55]{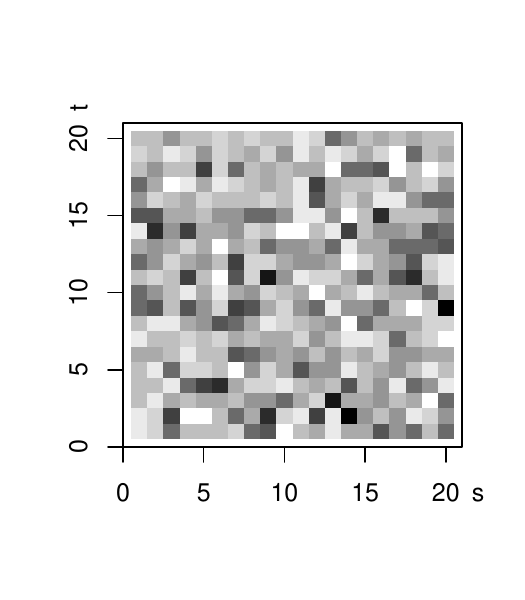}
\qquad
b)\includegraphics[viewport=0 45 335 240, clip=, scale=0.55]{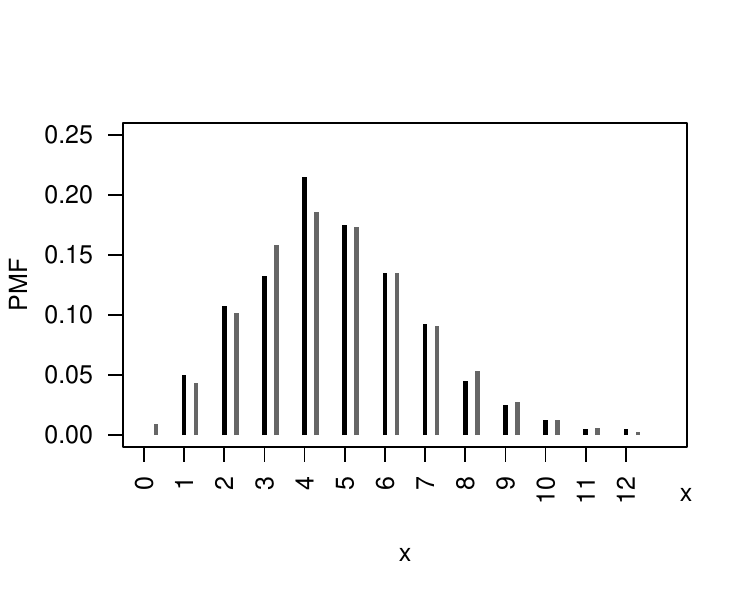}
\caption{Yeast count data discussed in Section~\ref{section: Numerical Illustrations}: a) gray-scale plot of data, where counts range between~0 (white) and~12 (black), and b) sample PMF (black) compared to $\poi(4.68)$-PMF (gray).}
\label{figure: yeast}
\end{figure}

As a first illustration of our novel INMA models for count random fields, let us briefly investigate a data set dating back to \citet[Table~I]{student06}, which was also analyzed by \citet{ghodsi12,ghodsi15}. The data being plotted in Figure~\ref{figure: yeast}~a) are counts of yeast cells over 1~mm\textsuperscript{2} divided into 400 squares (determined by using a haemacytometer), more precisely of yeast cells that were killed by adding a little mercuric chloride to the water \citep[p.~355]{student06}. The data exhibit only mild spatial dependence, with $\hat{\rho}(0,1)\approx 0.019$, $\hat{\rho}(1,0)\approx 0.034$, and $\hat{\rho}(1,1)\approx 0.039$. The marginal distribution, in turn, is similar to a Poisson distribution, see Figure~\ref{figure: yeast}~b), where the sample mean $\hat{\mu}_X=4.68$ is close to the sample variance $\hat{\sigma}_X^2\approx 4.47$, i.e., the yeast counts are nearly equidispersed with dispersion ratio $\hat{\sigma}_X^2/\hat{\mu}_X\approx 0.955$.

\smallskip
In \citet{ghodsi12,ghodsi15}, an INAR$(1,1)$ model with Poisson innovations (i.e., having the model equation $X_{s,t} = \alpha_{10}\circ X_{s-1,t} + \alpha_{01}\circ X_{s-1,t} + \alpha_{11}\circ X_{s-1,t} + \varepsilon_{s,t}$) was fitted to the data by the Yule--Walker (YW) approach, i.e., by computing the parameter estimates from sample mean and ACF. Since the estimate for~$\alpha_{11}$ was outside the admissible range $[0,1]$, the final model fit is given by $X_{s,t} = 0.040 \circ X_{s-1,t} + 0.021 \circ X_{s,t-1} + \varepsilon_{s,t}$ with $\varepsilon_{s,t}\sim\poi(4.567)$ (and with slightly deviating estimates in \citealp{ghodsi15}). But while the innovations are Poisson distributed, the INAR$(1,1)$'s observations are not --- in fact, their distribution is not known. Therefore, our Poisson INMA$(1,1)$--independence model, recall Examples~\ref{example: inma poisson} and~\ref{example: inma independence}, appears to be an attractive alternative. Estimating its model parameters again by a YW approach, i.e., numerically solving the three equations
    \[
        \hat{\rho}(k,l) = \sum_{(i,j)\in\mathcal{S}_{00}\cap\mathcal{S}_{kl}} \hat{\beta}_{i,j}\hat{\beta}_{i-k,j-l}/\hat{\beta}_\bullet \qquad \text{for } (k,l)\in \{(0,1), (1,0), (1,1)\}
    \]
in $\hat{\beta}_{0,1},\hat{\beta}_{1,0},\hat{\beta}_{1,1}$ with $\beta_{0,0}:=1$, we obtain the model fit
$$
X_{s,t} = \varepsilon_{s,t} + 0.020\circ \varepsilon_{s-1,t} + 0.036\circ \varepsilon_{s-1,t} + 0.043\circ \varepsilon_{s-1,t}
\quad\text{with } \varepsilon_{s,t}\sim\poi(4.257),
$$
which does not only exhibit a Poisson marginal distribution, but which also has another practical advantage compared to the INAR$(1,1)$ model: As discussed in Remark~\ref{remark: inma simulation}, the simulation of an INMA random field is computationally much cheaper than that of INAR random fields as burn-ins can be avoided. This, in turn, would be attractive if one wants to continue the above data analysis with, e.g., a parametric bootstrap approach.

\begin{remark}
Note that the YW approach for parameter estimation is easily implemented here due to the closed-form expression for the ACF of the specified INMA$(1,1)$ model (i.e., an independence model with Poisson innovations). In general, YW estimation requires a fully specified INMA model as the ACvF \eqref{eq: inma sacf} depends on the joint distribution of the counting series. This remains true even in the case of Poisson distributed innovations. It also transfers to the conditional moments according to Corollary~\ref{corollary: poisson inma dependence}, which might be used for developing a conditional least squares approach for parameter estimation. Maximum likelihood estimation, in turn, might become less feasible, because likelihood computation is quite demanding as INMA models are non-Markovian, recall the discussion around equations \eqref{eq: INMA representation}--\eqref{equation: pgf Y} as well as \eqref{equation: pgf Z independence}--\eqref{equation: pgf Y spread}. This shows the need for more elaborate methods for parameter estimation which should be addressed in future research.
\end{remark}

\smallskip
This leads us to a second illustration of our novel INMA models for count random fields. If one wants to simulate spatial data with an INAR-like dependence structure (e.g., for power analyses like in \citet{wei_kim_24}), but if one also wants to keep the required computing time at a feasible level (while knowing the actual marginal distribution), a possible solution is to use an ``INMA approximation'' of the INAR random field. More precisely, we speak of a $q$\textsuperscript{th}-order INMA approximation if the MA-parameters of an INMA$(q,q)$ model (with $\beta_{0,0}:=1$) are chosen such that it has the same ACF $\rho(k,l)$ with $0\leq k,l,\leq q$ as the aspired INAR model. Let us numerically illustrate our proposal for the case of an INAR$(1,1)$ random field with parameters $\alpha_{10} = \alpha_{01} = \alpha_{11}=0.2$. We refer to its ACF as the ``target ACF'' as we aim to find INMA parametrizations leading to the same ACF values for lags $0\leq k,l,\leq q$. Using the Poisson INMA$(q,q)$--independence model from Example~\ref{example: inma independence} for approximation, we compute the parameter values
\begin{itemize}
    \item $\beta_{0,1}=\beta_{1,0}\approx 0.383$, and $\beta_{1,1}\approx 0.783$ for $q=1$;
    
    \item $\beta_{0,1}=\beta_{1,0}\approx 0.342$, $\beta_{0,2}=\beta_{2,0}\approx 0.104$, $\beta_{1,1}\approx 0.623$, $\beta_{1,2}=\beta_{2,1}\approx 0.323$, and $\beta_{2,2}\approx 0.403$ for $q=2$;

    \item $\beta_{0,1}=\beta_{1,0}\approx 0.412$, $\beta_{0,2}=\beta_{2,0}\approx 0.032$, $\beta_{0,3}=\beta_{3,0}\approx 0.024$, $\beta_{1,1}\approx 0.809$, $\beta_{1,2}=\beta_{2,1}\approx 0.235$, $\beta_{1,3}=\beta_{3,1}\approx 0.110$, $\beta_{2,2}\approx 0.184$, $\beta_{2,3}=\beta_{3,2}\approx 0.159$, and $\beta_{3,3}\approx 0.188$ for $q=3$.
\end{itemize}
The resulting ACF values are summarized in Table~\ref{table: INMA approximation}, where we only show the lags $(k,l)$ with $0\leq k,|l|\leq 3$ as the missing ACF values for $k<0$ follow from the symmetry $\rho(-k,l)=\rho(k,-l)$. In this table, the approximation region is highlighted by gray background, and the ACF values being matched exactly by the definition of the INMA approximation are printed in bold font. It gets clear that we can successfully extend the INMA approximation to the required approximation quality. In particular, the gray values for $l<0$ are also quite close to those of the INAR random field, although the INMA approximation is defined with respect to the ACF values with $l\geq 0$ only. In the presented example, already the 2\textsuperscript{nd}-order but particularly the 3\textsuperscript{rd}-order approximation leads to a close approximation of the INAR$(1,1)$'s autocorrelation structure and might thus be sufficient for simulation purposes.

\begin{table}[t]
\centering\small
\caption{$q$\textsuperscript{th}-order INMA approximation of INAR$(1,1)$ random field, where lags with $k,|l|\leq q$ highlighted by gray background, and exactly matching ACF values $\rho(k,l)$ by bold font.}
\label{table: INMA approximation}

\smallskip
\begin{tabular}{@{}c@{\quad}c@{}}
\textbf{INAR$(1,1)$ target ACF} & \textbf{$q$\textsuperscript{th}-order INMA approximation}\\
& $q=1$\\[1ex]
\begin{tabular}{@{}r@{\qquad}cccc@{}}
\toprule
$l\ \setminus\ k$ & $0$ & $1$ & $2$ & $3$ \\[2ex]
$3$ & 0.019 & 0.044 & 0.057 & 0.046 \\
$2$ & 0.072 & 0.129 & 0.113 & 0.057 \\
$1$ & \highl\bf 0.268 & \highl\bf 0.307 & 0.129 & 0.044 \\
$0$ & \highl\bf 1.000 & \highl\bf 0.268 & 0.072 & 0.019 \\
$-1$ & \highl 0.268 & \highl 0.072 & 0.019 & 0.005 \\
$-2$ & 0.072 & 0.019 & 0.005 & 0.001 \\
$-3$ & 0.019 & 0.005 & 0.001 & 0.000 \\
\bottomrule
\end{tabular}
&
\begin{tabular}{@{}r@{\qquad}cccc@{}}
\toprule
$l\ \setminus\ k$ & $0$ & $1$ & $2$ & $3$ \\[2ex]
$3$ & 0.000 & 0.000 & 0.000 & 0.000 \\
$2$ & 0.000 & 0.000 & 0.000 & 0.000 \\
$1$ & \highl\bf 0.268 & \highl\bf 0.307 & 0.000 & 0.000 \\
$0$ & \highl\bf 1.000 & \highl\bf 0.268 & 0.000 & 0.000 \\
$-1$ & \highl 0.268 & \highl 0.058 & 0.000 & 0.000 \\
$-2$ & 0.000 & 0.000 & 0.000 & 0.000 \\
$-3$ & 0.000 & 0.000 & 0.000 & 0.000 \\
\bottomrule
\end{tabular}
\\
\\
& $q=2$\\[1ex]
\begin{tabular}{@{}r@{\qquad}cccc@{}}
\toprule
$l\ \setminus\ k$ & $0$ & $1$ & $2$ & $3$ \\[2ex]
$3$ & 0.019 & 0.044 & 0.057 & 0.046 \\
$2$ & \highl\bf 0.072 & \highl\bf 0.129 & \highl\bf 0.113 & 0.057 \\
$1$ & \highl\bf 0.268 & \highl\bf 0.307 & \highl\bf 0.129 & 0.044 \\
$0$ & \highl\bf 1.000 & \highl\bf 0.268 & \highl\bf 0.072 & 0.019 \\
$-1$ & \highl 0.268 & \highl 0.072 & \highl 0.019 & 0.005 \\
$-2$ & \highl 0.072 & \highl 0.019 & \highl 0.005 & 0.001 \\
$-3$ & 0.019 & 0.005 & 0.001 & 0.000 \\
\bottomrule
\end{tabular}
&
\begin{tabular}{@{}r@{\qquad}cccc@{}}
\toprule
$l\ \setminus\ k$ & $0$ & $1$ & $2$ & $3$ \\[2ex]
$3$ & 0.000 & 0.000 & 0.000 & 0.000 \\
$2$ & \highl\bf 0.072 & \highl\bf 0.129 & \highl\bf 0.113 & 0.000 \\
$1$ & \highl\bf 0.268 & \highl\bf 0.307 & \highl\bf 0.129 & 0.000 \\
$0$ & \highl\bf 1.000 & \highl\bf 0.268 & \highl\bf 0.072 & 0.000 \\
$-1$ & \highl 0.268 & \highl 0.098 & \highl 0.019 & 0.000 \\
$-2$ & \highl 0.072 & \highl 0.019 & \highl 0.003 & 0.000 \\
$-3$ & 0.000 & 0.000 & 0.000 & 0.000 \\
\bottomrule
\end{tabular}
\\
\\
& $q=3$\\[1ex]
\begin{tabular}{@{}r@{\qquad}cccc@{}}
\toprule
$l\ \setminus\ k$ & $0$ & $1$ & $2$ & $3$ \\[2ex]
$3$ & \highl\bf 0.019 & \highl\bf 0.044 & \highl\bf 0.057 & \highl\bf 0.046 \\
$2$ & \highl\bf 0.072 & \highl\bf 0.129 & \highl\bf 0.113 & \highl\bf 0.057 \\
$1$ & \highl\bf 0.268 & \highl\bf 0.307 & \highl\bf 0.129 & \highl\bf 0.044 \\
$0$ & \highl\bf 1.000 & \highl\bf 0.268 & \highl\bf 0.072 & \highl\bf 0.019 \\
$-1$ & \highl 0.268 & \highl 0.072 & \highl 0.019 & \highl 0.005 \\
$-2$ & \highl 0.072 & \highl 0.019 & \highl 0.005 & \highl 0.001 \\
$-3$ & \highl 0.019 & \highl 0.005 & \highl 0.001 & \highl 0.000 \\
\bottomrule
\end{tabular}
&
\begin{tabular}{@{}r@{\qquad}cccc@{}}
\toprule
$l\ \setminus\ k$ & $0$ & $1$ & $2$ & $3$ \\[2ex]
$3$ & \highl\bf 0.019 & \highl\bf 0.044 & \highl\bf 0.057 & \highl\bf 0.046 \\
$2$ & \highl\bf 0.072 & \highl\bf 0.129 & \highl\bf 0.113 & \highl\bf 0.057 \\
$1$ & \highl\bf 0.268 & \highl\bf 0.307 & \highl\bf 0.129 & \highl\bf 0.044 \\
$0$ & \highl\bf 1.000 & \highl\bf 0.268 & \highl\bf 0.072 & \highl\bf 0.019 \\
$-1$ & \highl 0.268 & \highl 0.086 & \highl 0.021 & \highl 0.004 \\
$-2$ & \highl 0.072 & \highl 0.021 & \highl 0.006 & \highl 0.001 \\
$-3$ & \highl 0.019 & \highl 0.004 & \highl 0.001 & \highl 0.000 \\
\bottomrule
\end{tabular}
\end{tabular}
\end{table}

\bigskip
As a final illustration, let us briefly return to the multilateral INMA model \eqref{eq: model inma multilateral}, the ACF of which is calculated in complete analogy to that of the unilateral model, recall Remark~\ref{remark: multilateral}. In what follows, we always focus on Poisson INMA--independence models according to Example~\ref{example: inma independence} for simplicity. Table~\ref{table: multilateral INMA} compares the ACF values of the unilateral INMA$(1,1)$ model (left block) to the multilateral INMA$(1,1)$ model (central block), where we always use the same parametrization, namely $\beta_{0,0}=1$ and~$\beta_{i,j}=0.2$ otherwise. It gets clear that the multilateral model has both a longer memory and considerably stronger ACF values, despite using analogous parametrizations. Both observations are plausible in view of the more sophisticated dependence structure of the multilateral INMA model, recall Figure~\ref{fig: multilateral structure}. Since the multilateral INMA$(1,1)$ model has a memory of length~2, we can apply a 2\textsuperscript{nd}-order INMA approximation as described before, i.e., we determine the parameters of a unilateral INMA$(2,2)$ model such that it has the same ACF values for lags $0\leq k,l\leq 2$. The corresponding parameter values are $\beta_{0,1}=\beta_{1,0}\approx 0.326$, $\beta_{0,2}=\beta_{2,0}\approx 0.086$, $\beta_{1,1}\approx 0.383$, $\beta_{1,2}=\beta_{2,1}\approx 0.061$, and $\beta_{2,2}\approx 0.036$, and the resulting ACF values are shown in the right block of Table~\ref{table: multilateral INMA}. But despite having the same ACF values for lags $0\leq k,l\leq 2$, the multilateral INMA$(1,1)$ model still has clearly larger ACF values for $l<0$. So although having the same number of model parameters, the multilateral INMA$(1,1)$ model leads to an overally stronger dependence structure, which motivates to analyze this model in more detail in future research.  

\begin{table}[t]
\centering\small
\caption{ACF values $\rho(k,l)$ of some unilateral and multilateral INMA models, see Section~\ref{section: Numerical Illustrations}.}
\label{table: multilateral INMA}

\smallskip
\begin{tabular}{@{}rcccccccccccc@{}}
\toprule
&& \multicolumn{3}{c}{\bf unilat.\ INMA$(1,1)$} && \multicolumn{3}{c}{\bf multilat.\ INMA$(1,1)$} && \multicolumn{3}{c}{\bf unilat.\ INMA$(2,2)$} \\
$l\ \setminus\ k$ && $0$ & $1$ & $2$ && $0$ & $1$ & $2$ && $0$ & $1$ & $2$ \\
\cmidrule{1-1}\cmidrule{3-5}\cmidrule{7-9}\cmidrule{11-13}
$2$ && 0.000 & 0.000 & 0.000 && 0.046 & 0.031 & 0.015 && 0.046 & 0.031 & 0.015 \\
$1$ && 0.150 & 0.125 & 0.000 && 0.215 & 0.185 & 0.031 && 0.215 & 0.185 & 0.031 \\
$0$ && 1.000 & 0.150 & 0.000 && 1.000 & 0.215 & 0.046 && 1.000 & 0.215 & 0.046 \\
$-1$ && 0.150 & 0.025 & 0.000 && 0.215 & 0.185 & 0.031 && 0.215 & 0.074 & 0.014 \\
$-2$ && 0.000 & 0.000 & 0.000 && 0.046 & 0.031 & 0.015 && 0.046 & 0.014 & 0.003 \\
\bottomrule
\end{tabular}
\end{table}

%%%%%%%%%%%%%%%%%%%%%%%%
\section{Conclusions and Future Research}
\label{section: Conclusions and Future Research}

In this article, we introduced an INMA-type random field model for counts being observed on a two-dimensional regular grid. While mainly focusing on a unilateral model structure, we also briefly considered the multilateral case. We derived closed-form expressions for both the INMA's marginal distribution and spatial dependence structure under arbitrary model order. In particular, general expressions for bivariate distributions (pgf) and autocovariances at a specified spatial lag $(k,l)$ are provided. We showed that the INMA random field can be equipped with a Poisson marginal distribution, which differs from existing INAR-type random-field models in the literature. We also illustrated that different and well-interpretable dependence structures are possible, by developing the special INMA cases of the independence and spread model, and by also comparing to a multilateral INMA model. The practical relevance of our novel INMA model for count random fields was demonstrated by a real-world data example on the one hand, and by the INMA approximation of an INAR random field's ACF on the other hand. The latter application was motivated by the fact that an INMA random field is simulated with considerably less computational effort than an INAR random field. Moreover, we compared the unilateral and multilateral INMA models with regard to their ACF. Despite using an analogous parametrization, we showed that the multilateral model has both a longer and an intensified memory, which is not observed in the same way for corresponding unilateral models.

\smallskip
While the present article focused on specifications and stochastic properties of INMA random fields, future research should be directed towards statistical inference. In particular, methods for parameter estimation and model diagnostics should be developed and investigated. Furthermore, the multilateral INMA model deserves further attention. In addition, it would be of practical relevance if also INAR-type random fields with a Poisson (or another common) marginal distribution could be developed. Moreover, the extension of INARMA approaches to random fields with other types of discrete-valued range (such as signed integers or bounded counts) would be an interesting task for future research. 
Finally, considering integer-valued random fields in the generalized linear model framework is an interesting direction for future research, paving the way for \emph{multivariate} count random field models in particular. For example, extensions of the count network autoregression models proposed by \citet{armillotta24,guo24a} and \citet{guo24b} to random fields are conceivable.

%%%%%%%%%%%%%%%%%%%%%%
\subsubsection*{Acknowledgments}
The authors thank the associate editor and the referees for their useful comments on an earlier draft of this article.

%%%%%%%%%%%%%%%%%%%%%%

\bibliographystyle{plainnat}
\bibliography{references}
%\vspace{12pt}

% \newpage
\appendix
\numberwithin{equation}{section}
%\numberwithin{table}{section}

\section{Derivations}
\label{Derivations}

\subsection{Proof of Proposition~\ref{prop: inma marginal properties}}
\label{app: Proof of Proposition prop: inma marginal properties}
The first expression is an immediate consequence of the fact that for all constants $\alpha\in[0,1]$ and random variables $X$, it holds $\alpha \circ X | X \sim \bin(X, \alpha)$ and hence, by the law of total expectation, $\bbe(\alpha\circ X) = \bbe(\bbe(\alpha \circ X| X)) = \alpha\, \bbe (X)$. Since $(\varepsilon_{s,t})$ is i.i.d., and since the thinnings applied to different~$\varepsilon_{s,t}$ are independent of each other as well, Lemma~1 of \citet{sil_wei_26} implies
\[
    \var(X_{s,t}) = \sum^{q_1}_{i=0} \sum^{q_2}_{j=0} \var(\beta_{i,j}\circ \varepsilon_{s-i,t-j}) = \sum^{q_1}_{i=0} \sum^{q_2}_{j=0} \Bbrackets{\beta_{i,j}^2 \sigma^2_\varepsilon + \beta_{i,j}(1-\beta_{i,j}) \mu_\varepsilon}, 
\]
which leads to the desired expression. Note that for a random variable $X$ and a constant $\alpha \in [0,1]$, it holds
\[
    \var(\alpha \circ X) = \var(\bbe(\alpha \circ X|X)) + \bbe(\var(\alpha \circ X|X)) = \var(\alpha \cdot X) + \alpha (1-\alpha) \bbe (X) \neq \var(\alpha \cdot X),
\]
see \citet[p.~17]{weiss18}, for instance.
With regard to the pgf, first we consider $\pgf_{\beta\circ\varepsilon}(u) = \bbe(u^{\beta\circ\varepsilon})$. Let $(Z_r) = (Z_r)_{r\in\bbn}$ denote the i.i.d.\ counting series corresponding to the binomial thinning $\beta \circ \varepsilon := \sum^\varepsilon_{r=1} Z_r$. Then,
\begin{align}
    \pgf_{\beta\circ\varepsilon}(u) &=  \bbe\brackets{\bbe(u^{\sum^\varepsilon_{r=1}Z_r}\, |\, \varepsilon)} 
    = \bbe\Big(\big(\bbe(u^{Z_1})\big)^\varepsilon\Big) \nonumber \\
    &= \bbe\Big(\big((1-\beta) u^0 + \beta\, u\big)^\varepsilon\Big)= \pgf_\varepsilon\big(1+\beta(u-1)\big). \label{eq: pgf thinning}
\end{align}
Hence, it follows that
\begin{align*}
    \pgf_X(u) = \bbe\brackets{\prod^{q_1}_{i=0}\prod^{q_2}_{j=0} u^{\beta_{i,j}\circ \varepsilon_{s-i,t-j}}} 
    &= \prod^{q_1}_{i=0}\prod^{q_2}_{j=0} \bbe\brackets{\bbe\brackets{u^{\beta_{i,j}\circ \varepsilon_{s-i,t-j}}\, \big|\, \varepsilon_{s-i,t-j}}} \\
    &= \prod^{q_1}_{i=0}\prod^{q_2}_{j=0} \pgf_\varepsilon\big(1+\beta_{i,j}(u-1)\big),
\end{align*}
and the proof of Proposition~\ref{prop: inma marginal properties} is complete.

\subsection{Proof of Theorem~\ref{theorem: inma dependence}}
\label{app: Proof of Theorem theorem: inma dependence}
Part~(i):\quad
Recall the sets $\mathcal{S}_{kl}\, :=\, \{(i+k,j+l) : 0 \leq i \leq q_1, 0 \leq j \leq q_2\}$ for $k,l\in\bbz$. 
By considering the three regions in Figure~\ref{fig: proof inma bivariate pgf}, we have
\begin{align*}
    \pgf(u_1,u_2;k,l) &= \bbe\big(u_1^{X_{s,t}} u_2^{X_{s-k,t-l}}\big) \\
    &= \prod_{(i,j)\in\mathcal{S}_{00}\setminus(\mathcal{S}_{00}\cap\mathcal{S}_{kl})} \bbe\big(u_1^{\beta_{i,j}\circ\varepsilon_{s-i,t-j}}\big)
    \cdot \prod_{(i,j)\in\mathcal{S}_{kl}\setminus(\mathcal{S}_{00}\cap\mathcal{S}_{kl})} \bbe\big(u_2^{\beta_{i-k,j-l}\circ\varepsilon_{s-i,t-j}}\big)  \\
    &\qquad \cdot\ \prod_{(i,j)\in\mathcal{S}_{00}\cap\mathcal{S}_{kl}} \bbe\big(u_1^{\beta_{i,j}\circ\varepsilon_{s-i,t-j}} u_2^{\beta_{i-k,j-l}\circ\varepsilon_{s-i,t-j}}\big).
\end{align*}
The first two expected values can be determined by using \eqref{eq: pgf thinning}. Concerning the last term appearing in the expression above (corresponding to the gray region in Figure~\ref{fig: proof inma bivariate pgf}), in order to shorten notations, we write $\varepsilon = \varepsilon_{s-i,t-j}$, $Z_r^{(i,j)}=Z_{s-i,t-j;r}^{(i,j)}$ and $Z_r^{(i-k,j-l)}=Z_{s-i,t-j;r}^{(i-k,j-l)}$, respectively. Then, we have 
\begin{align*}
    \bbe\big(u_1^{\beta_{i,j} \circ \varepsilon} u_2^{\beta_{i-k,j-l} \circ \varepsilon}\big) 
    &= \bbe\brackets{\bbe\brackets{u_1^{\sum^\varepsilon_{r=1}Z_r^{(i,j)}} u_2^{\sum^\varepsilon_{r=1} Z_r^{(i-k,j-l)}}\, \Big|\, \varepsilon}} \\
    &= \bbe\brackets{\brackets{\bbe(u_1^{Z_1^{(i,j)}} u_2^{Z_1^{(i-k,j-l)}})}^\varepsilon} \\
    &= \pgf_\varepsilon\Big(\bbe\big(u_1^{Z_1^{(i,j)}} u_2^{Z_1^{(i-k,j-l)}}\big)\Big)
\end{align*}
due to $(\varepsilon_{s,t})$ being i.i.d.
Since 
\begin{align*}
    \bbp(Z^{(i,j)}_{1}=1, Z^{(i-k,j-l)}_{1}=0) &= \beta_{i,j} -  \bbp(Z^{(i,j)}_{1}=Z^{(i-k,j-l)}_{1}=1), \\
    \bbp(Z^{(i,j)}_{1}=0, Z^{(i-k,j-l)}_{1}=1) &= \beta_{i-k,j-l} - \bbp(Z^{(i,j)}_{1}=Z^{(i-k,j-l)}_{1}=1), \\
    \bbp(Z^{(i,j)}_{1}=Z^{(i-k,j-l)}_{1}=0) &= 1 - \beta_{i,j} - \beta_{i-k,j-l} + \bbp(Z^{(i,j)}_{1}=Z^{(i-k,j-l)}_{1}=1),
\end{align*}
it holds that
\begin{align*}
    \bbe\big(u_1^{Z^{(i,j)}} u_2^{Z^{(i-k,j-l)}}\big) &= \sum^1_{a,b=0} u_1^a u_2^b\, \bbp(Z^{(i,j)}_{1}=a, Z^{(i-k,j-l)}_{1}=b) \\
    &= 1 + \beta_{i,j} (u_1-1) + \beta_{i-k,j-l} (u_2-1) \\
    &\qquad + \bbp(Z^{(i,j)}_{1}=Z^{(i-k,j-l)}_{1}=1)\, (u_1-1)(u_2-1).
\end{align*}
Putting everything together, the proof of Theorem~\ref{theorem: inma dependence}~(i) is complete.

\smallskip
Part~(ii):\quad
First of all, it holds that
\begin{align*}
    \gamma(k,l) &= \cov(X_{s,t},X_{s-k,t-l}) \\
    &= \sum^{q_1}_{i_1=0} \sum^{q_2}_{j_1=0} \sum^{q_1}_{i_2=0} \sum^{q_2}_{j_2=0} \cov(\beta_{i_1,j_1}\circ_{s,t}\varepsilon_{s-i_1,t-j_1},\ \beta_{i_2,j_2} \circ_{s-k,t-l}\varepsilon_{s-k-i_2,t-l-j_2}) \\
    &= \sum_{(i,j)\in\mathcal{S}_{00}\cap\mathcal{S}_{kl}} \cov(\beta_{i,j}\circ_{s,t}\varepsilon_{s-i,t-j},\ \beta_{i-k,j-l}\circ_{s-k,t-l}\varepsilon_{s-i,t-j}),
\end{align*}
since all terms with non-coinciding $\varepsilon$-s vanish, i.e., only the terms with $i_1=i_2+k$ and $j_1=j_2+l$ remain. Here, we set $i_1=i$ and $j_1=j$ (such that $i_2=i-k$ and $j_2=j-l$) to obtain the third equation.
To simplify notation, we write $\varepsilon:=\varepsilon_{s-i,t-j}$ in the remainder of this proof, as well as $(Z^{(i,j)}_{r_1})$ and $( Z^{(i-k,j-l)}_{r_2})$ for the corresponding counting series if thinnings are applied to~$\varepsilon$. Considering the remaining covariances, by the law of total covariance and the given independence properties, it follows that
\begin{align*}
    &\cov(\beta_{i,j}\circ_{s,t} \varepsilon,\, \beta_{i-k,j-l}\circ_{s-k,t-l}\varepsilon) \\
    &\qquad = \cov\Big(\bbe(\beta_{i,j}\circ_{s,t} \varepsilon\, |\, \varepsilon), \bbe(\beta_{i-k,j-l} \circ_{s-k,t-l}\varepsilon\, |\, \varepsilon)\Big) \\
    &\qquad \qquad \qquad + \bbe\Big(\cov(\beta_{i,j}\circ_{s,t}\varepsilon,\, \beta_{i-k,j-l} \circ_{s-k,t-l}\varepsilon\, |\, \varepsilon)\Big) \\
    &\qquad = \cov(\beta_{i,j}\, \varepsilon, \beta_{i-k,j-l}\, \varepsilon) + \bbe\biggl(\cov\Bigl(\sum^\varepsilon_{r_1=1} Z^{(i,j)}_{r_1}, \sum^\varepsilon_{r_2=1} Z^{(i-k,j-l)}_{r_2}\, \Big|\, \varepsilon\Bigr)\biggr) \\
    &\qquad = \beta_{i,j}\beta_{i-k,j-l} \,\sigma^2_\varepsilon + \bbe\Bigl(\varepsilon\cdot \cov\bigl(Z^{(i,j)}_{1},\, Z^{(i-k,j-l)}_{1}\bigr)\Bigr) \\
    &\qquad = \beta_{i,j}\beta_{i-k,j-l} \,\sigma^2_\varepsilon +\mu_\varepsilon\, \Bigl(\bbp\bigl(Z^{(i,j)}_{1}= Z^{(i-k,j-l)}_{1}=1\bigr) - \bbp\bigl(Z^{(i,j)}_{1}=1\bigr) \bbp\bigl(Z^{(i-k,j-l)}_{1}=1\bigr)\Bigr) \\
    &\qquad = \beta_{i,j}\beta_{i-k,j-l} \,\sigma^2_\varepsilon +\mu_\varepsilon\, \Bigl(\bbp\bigl(Z^{(i,j)}_{1}= Z^{(i-k,j-l)}_{1}=1\bigr) - \beta_{i,j} \beta_{i-k,j-l}\Bigr).
\end{align*}
This completes the proof of Theorem~\ref{theorem: inma dependence}.

\subsection{Proof of Corollary~\ref{corollary: poisson inma dependence}}
\label{app: Proof of Corollary corollary: poisson inma dependence}
The ACvF according to Theorem~\ref{theorem: inma dependence}~(ii) simplifies considerably in the Poisson case, due to the equidispersion property $\sigma^2_\varepsilon = \mu_\varepsilon$:
\begin{equation*}
    \gamma(k,l) = \mu_\varepsilon \sum_{(i,j)\in\mathcal{S}_{00}\cap\mathcal{S}_{kl}} \bbp(Z_{s-i,t-j;1}^{(i,j)}=Z_{s-i,t-j;1}^{(i-k,j-l)}=1).
\end{equation*}
Furthermore, $\gamma(0,0) = \mu_\varepsilon \beta_\bullet$ according to Example~\ref{example: inma poisson}, which implies the expression for the ACF in part~(i).

\smallskip
Concerning the bivariate pgf in part~(ii), \eqref{eq: pgf poisson} and \eqref{PoiINMA_acf} imply that
    \begin{align*}
        \pgf(u_1,u_2;k,l) = & \prod_{(i,j)\in\mathcal{S}_{00}\setminus(\mathcal{S}_{00}\cap\mathcal{S}_{kl})} \exp\big(\mu_\varepsilon\beta_{i,j}(u_1-1)\big) \\
        &\qquad\cdot \prod_{(i,j)\in\mathcal{S}_{kl}\setminus(\mathcal{S}_{00}\cap\mathcal{S}_{kl})} \exp\big(\mu_\varepsilon\beta_{i-k,j-l}(u_2-1)\big) \\
        &\qquad \cdot \prod_{(i,j)\in\mathcal{S}_{00}\cap\mathcal{S}_{kl}} \exp\Big[\mu_\varepsilon\Big(\beta_{i,j}(u_1-1) + \beta_{i-k,j-l}(u_2-1) \\
        &\hspace{30mm} + \bbp(Z^{(i,j)}_{s,t;1}=Z^{(i-k,j-l)}_{s,t;1}=1) \cdot(u_1-1)(u_2-1)\Big)\Big] \\
        &= \exp\big(\mu_\varepsilon\beta_\bullet(u_1-1)\big) \cdot \exp\big(\mu_\varepsilon\beta_\bullet(u_2-1)\big) \\
        &\qquad \cdot \exp\Big(\mu_\varepsilon(u_1-1)(u_2-1) \sum_{(i,j)\in\mathcal{S}_{00}\cap\mathcal{S}_{kl}} \bbp(Z^{(i,j)}_{s,t;1}=Z^{(i-k,j-l)}_{s,t;1}=1)\Big) \\
        &= \exp\big(\mu_X (u_1+u_2-2)\big) \cdot \exp\big(\mu_X(u_1-1)(u_2-1)\rho(k,l)\big).
    \end{align*}
The proof of the conditional moments in Corollary~\ref{corollary: poisson inma dependence}~(iii) goes along the same lines as that of Theorem~3.1 in \citet{weiss08}, where $\rho(k)$ is replaced by $\rho(k,l)$, and where the definition of $\beta_\bullet$ is adjusted accordingly, because aside from the aforementioned terms, the pgfs under consideration do not differ. With the same changes, also the proof of part~(iv) corresponds to that of \citet[Theorem~3.3]{weiss08}. So the proof of Corollary~\ref{corollary: poisson inma dependence} is complete.

\end{document}